\documentstyle[11pt]{article}
\newfont{\Bbb}{msbm10 scaled\magstephalf}
\begin{document}
\title{Compact Composition Operators on the Bloch Space
in Bounded Symmetric Domains}
\author{Zehua Zhou\hspace*{4mm}Yan Liu\\
Department of Mathematics, Tianjin University, Tianjin 300072\\
E-mail: zehuazhou2003@yahoo.com.cn}
\date{}
\maketitle \footnote{Supported in part by the National Natural
Science Foundation of China (Grand No.10371091), and LiuHui Center
for Applied Mathematics, Nankai University \& Tianjin University}

\begin{abstract} Let $\Omega$ be a
bounded symmetric domain except the two exceptional domains of
$\mbox{\Bbb C}^N$ and $\phi$ a holomorphic self-map of $\Omega.$ This paper gives a sufficient
and necessary condition for the composition operator $C_{\phi}$
induced by $\phi$ to be compact on the Bloch space $\beta(\Omega)$.
\end{abstract}

\hspace{4mm}{\bf Key words}\hspace{2mm}Bloch space,
Bounded symmetric domains, Composition operator, Bergman metric.

\hspace{4mm}{\bf 2000 Mathematics Subject Classification}
\hspace{2mm}47B38\hspace{2mm}47B33\hspace{2mm}32A30\hspace{2mm}32A37

\section{Introduction}
\newtheorem{Definition}{Definition}
\newtheorem{Lemma}{Lemma}
\newtheorem{Theorem}{Theorem}
\newtheorem{Proposition}{Proposition}

Let $\cal D$ be a bounded homogeneous domain in
$\mbox{\Bbb C}^{N}.$ The class of all holomorphic functions with domain
$\cal D$ will be denoted by $H(\cal D).$ Let $\phi $ be a holomorphic
self-map of $\cal D$. For $f\in H(\cal D),$
we denote the composition $f\circ\phi$ by $C_{\phi}f$ and call $C_{\phi}$
the composition operator induced by $\phi$.

Let $K(z,z)$ be the Bergman kernel function of $\cal D$. The Bergman metric
$H_{z}(u,u)$ in $\cal D$ is defined by
$$H_{z}(u,u)=\displaystyle\frac{1}{2}
\sum\limits^{N}_{l,k=1}
\displaystyle\frac{\partial^{2}\log K(z,z)}{\partial z_l
\partial {\overline {z}}_{k}}u_l{\overline u}_{k},$$
where $z\in \cal D$
and $u=(u_{1},\ldots,u_{N})\in \mbox{\Bbb C}^{N}.$

Following Timoney [1], we say that $f\in H(\cal D)$ is in the Bloch space
$\beta(\cal D),$ if \begin{equation}
\|f\|_{\beta(\cal D)}=\sup_{z\in \cal D}Q_{f}(z)<\infty,\label{a}
\end{equation}
where$$Q_{f}(z)=\sup\left\{\displaystyle\frac
{|\bigtriangledown f(z)u|}{H^{\frac{1}{2}}_{z}(u,u)}:
u\in \mbox{\Bbb C}^{N}-\{0\}\right\},$$ and $\bigtriangledown f(z)
=\left(\frac{\partial f(z)}{\partial z_{1}}, \ldots,
\frac{\partial f(z)}{\partial z_{N}} \right), \bigtriangledown f(z)u
=\sum\limits^{N}_{l=1}\frac{\partial f(z)}{\partial z_l}u_l.$

Let $D$ be the unit disk in $\mbox{\Bbb C}.$
Madigan and Matheson [2] proved that $C_{\phi}$ is always bounded
on $\beta(D).$ They also gave the
sufficient and necessary condition that $C_{\phi}$ is compact on
$\beta(D).$

Recently, Shi and Luo [3] proved that $C_{\phi}$ is always bounded
on $\beta(\cal D),$ where $\cal D$ is a bounded homogeneous domain
in $\mbox{\Bbb C}^{N}.$ They also gave a sufficient condition for
$C_{\phi}$ to be compact on $\beta(\cal D)$ (i.e., Lemma 3). So
this result leads us to ask whether the condition is also
necessary. Zhou and Shi [4] give an affirmative answer to this
question for classical bounded symmetric domains. In fact, the
original purpose of the typescript for [4], is to answer the
question in the bounded symmetric domains of $\mbox{\Bbb C}^N$,
but in the final proof, the referee point out a fatal mistake
which we aren't able to correct that time, upon the suggestion of
the referee, we cancel out the last part and published in the form
of [4] which discussed only in the four types classical bounded
symmetric domains. The following paper overcome difficulty and
solve the problem, for the method used here we essentially follow
[4], but some new techniques have been used.

By using Cartan's list, all irreducible bounded symmetric domains
are divided into six types. The first four types of irreducible
domains are called the classical bounded symmetric domains. The
other two types, called exceptional domains, consist of one domain
each (a 16- and a 27-dimensional domain).

In what follows, $\Omega$ denotes a bounded symmetric domain
except the two exceptional domains of $\mbox{\Bbb C}^N$, and
$\phi$ a holomorphic self-map of $\Omega.$ If $U=(u_{kl})_{m\times
n}$ is a $m\times n$ complex matrix, write $u=(u_{11}, \ldots,
u_{1n}, \ldots, u_{m1}, \ldots, u_{mn})$ as the corresponding
vector of matrix $U$ and $\overline u'$ is the conjugate transpose
of $u.$ $C$ is a positive constant, not necessarily the same at
each occurrence.

In this paper, we will give a sufficient and
necessary condition that the composition operator $C_{\phi}$ is compact
on $\beta(\Omega).$

Let $A=(a_{jk})_{m\times n}, B=(b_{lr})_{p\times q}$. The Kronecker
product of $A$ and $B,$ defined by
$A\times B=C=(c_{jlkr}),$ is a $mp\times nq$ matrix,
where $c_{jlkr}=a_{jk}b_{lr}$.

It is well known [5] that the classical bounded symmetric domains
$R_I, R_{II}, R_{III},$ and $R_{IV}$ can be expressed as follows:

$R_{I}(m,n)=\{Z: Z\mbox{ is a } m\times n
\mbox{ complex matrix, }
I_{m}-Z{\overline Z}'>0\},$ where $I_{m}$ is the $m\times m$ identity
matrix ($m\leq n$).

$R_{II}(p)=\{Z: Z\mbox{ is a } p\times p\mbox{ symmetric matrix } Z=Z',
I_{p}-Z{\overline Z}>0\}.$

$R_{III}(q)=\{Z: Z\mbox{ is a } q\times q\mbox{ antisymmetric matrix }Z=-Z,'
I_{q}+Z{\overline Z}>0\}.$

$R_{IV}(N)=\{z: z=(z_{1}, \ldots, z_{N}),
1+|zz'|^{2}-2z{\overline z}'>0, |zz'|<1\}.$

Their Bergman metrics are the following respectively [6]:
\begin{equation}H^{I}_{z}(u,u)=(m+n)u(I_{m}-Z{\overline Z}')^{-1}\times
(I_{n}-{\overline Z}'Z)^{-1}{\overline u}',\label{b}\end{equation}
where $Z\in R_{I}(m,n)$ and
$U$ is a $m\times n$ complex matrix, $u$ is the corresponding vector of
$U,$  $\overline u'$ is the conjugate transpose of $u.$
\begin{equation}H^{II}_{z}(u,u)=(p+1)u(I_{p}-Z{\overline Z})^{-1}\times
(I_{p}-{\overline Z}Z)^{-1}{\overline u}',\label{c}\end{equation}
where $Z\in R_{II}(p)$ and
$U$ is a $p\times p$ symmetric complex matrix,
$u$ is the corresponding vector of $U.$
\begin{equation}H^{III}_{z}(u,u)=2(q-1)u(I_{q}+Z{\overline Z})^{-1}\times
(I_{q}+{\overline Z}Z)^{-1}{\overline u}',\label{d}\end{equation}
where $Z\in R_{III}(q)$ and $U$ is a $q\times q$ anti-symmetric
complex matrix, $u$ is the corresponding vector of $U.$
\begin{eqnarray}H^{IV}_{z}(u,u)&=&\displaystyle\frac{2N}
{(1+|zz'|^{2}-2z{\overline z}')^{2}}
u\Biggl[(1+|zz'|^{2}-2z{\overline z}')I_N\nonumber\\
&&-2\left(\begin{array}{c}z \\
{\overline z}\end{array}\right)'\left(\begin{array}{lr}1-2|z|^{2}&
\overline{zz'}\\zz'& -1\end{array}\right)
\overline{\left(\begin{array}{c}z\\
{\overline z}\end{array}\right)}\Biggr]
{\overline u}'.\label{e}\end{eqnarray}
where $z\in R_{IV}(N)$ and $u\in\mbox{\Bbb C}^N.$

Our main result is the following:

{\bf Theorem}\hspace{2mm} Let $\Omega\subset\mbox{\Bbb C}^N$ be a
bounded symmetric domain except the two exceptional domains and $\phi$
a holomorphic self-map of $\Omega.$ Then $C_{\phi}$ is compact
on the Bloch space $\beta(\Omega)$ if and only if for every
$\varepsilon>0$, there exists a $\delta >0,$ such that \begin{equation}
\displaystyle\frac{H_{\phi(z)}\left(J\phi(z)u, J\phi(z)u\right)}
{H_{z}(u,u)}<\varepsilon,\label{f}\end{equation}for all $u\in \mbox{\Bbb C}^{N}-\{0\}$
whenever $dist(\phi(z), \partial\Omega)<\delta,$
where $H_z(u, u)$ is the
Bergman metric of $\Omega.$

{\bf Remark}\hspace{2mm} It is well known that the unit ball and
unit polydiscs are bounded symmetric domains, so the above result
holds in the unit balls and unit polydiscs. Furthermore, we can
also obtain Theorem 2 in [2].

\section{Some Lemmas}
In order to prove the Theorem, we need the following Lemmas.
\begin{Lemma}([1], Theorem 2.12)
Let ${\cal D}\subset\mbox{\Bbb C}^N$ be a bounded homogeneous domain. Then
there exists a constant $C$ depending only on
$\cal D$, such that $H_{\phi(z)}\left(J\phi(z)u, J\phi(z)u\right)
\leq CH_{z}(u,u)$, for each $z\in \cal D,$ whenever $\phi$ holomorphically
maps $\cal D$ into itself. Here $H_{z}(u,u)$  denotes the Bergman metric on
$\cal D,$ $J\phi(z)=\left(
\displaystyle\frac{\partial \phi_{l}(z)}
{\partial z_{k}}\right)_{1\leq l,k\leq N}$ denotes the Jacobian matrix of
$\phi ,$ and $J\phi(z)u$ denotes a vector, whose $l$th component is
$\left(J\phi(z)u\right)_{l}=\sum\limits^{N}_{k=1}
\displaystyle\frac{\partial \phi_{l}(z)}
{\partial z_{k}}u_{k},$ $l=1,2,\ldots,N.$\end{Lemma}

\begin{Lemma} ([3], Lemma 3) Let $\cal D$ be a bounded homogeneous
domain in $\mbox{\Bbb C}^{N}.$ Then $C_{\phi}$ is compact
on $\beta(\cal D)$ if and only if for any bounded sequence $\{f_{k}\}$
in $\beta(\cal D)$ which
converges to zero uniformly on compact subsets of $\cal D$, we have
$\|f_{k}\circ \phi\|_{\beta(\cal D)}\to 0,$ as $k \to\infty.$\end{Lemma}

\begin{Lemma} ([3], Theorem 3) If $\phi: \cal D\to\cal D$
is a holomorphic self-map, where $\cal D$ is a bounded homogeneous
domain in $\mbox{\Bbb C}^{N}.$
Then $C_{\phi}$ is compact on $\beta(\cal D)$ if for every $\varepsilon>0$,
there exists a $\delta >0,$ such that $$
\displaystyle\frac{H_{\phi(z)}\left(J\phi(z)u, J\phi(z)u\right)}
{H_{z}(u,u)}<\varepsilon,$$for all $u\in \mbox{\Bbb C}^{N}-\{0\}$
whenever $dist(\phi(z), \partial\cal D)<\delta.$\end{Lemma}

\begin{Lemma} ([4], Lemma 4) Let $\cal D$ be a bounded homogeneous domain of
$\mbox{\Bbb C}^N,$ and let $T(z,z)$ denote its metric matrix. If
$T(0,0)=\lambda I_N,$ where $\lambda$ is a constant depending only
on $\cal D,$ then a holomorphic function $f$ on $\cal D$ is in
$\beta(\cal D)$ if and only if
\begin{equation}\sup\limits_{z\in\cal D}\left\{\bigtriangledown f(z)T^{-1}(z,z)
\overline{\bigtriangledown f(z)}'\right\}<\infty.\label{g}\end{equation}
If (\ref{g}) holds, then there exists a constant $C$ depending only on $\cal D$
such that $$\|f\|_{\beta(\cal D)}\leq C\sup\limits_{z\in\cal D}
\left\{\bigtriangledown f(z)T^{-1}(z,z)
\overline{\bigtriangledown f(z)}'\right\}.$$\end{Lemma}

\begin{Lemma}([1], Proposition 4.5) Let $\cal D$ be a bounded
homogeneous domain in $\mbox{\Bbb C}^{N}.$
If $f$ is a bounded holomorphic function in $\cal D,$
then $f\in\beta(\cal D)$ and there exists a
constant $C$ depending only on $\cal D$, such that \
$$\|f\|_{\beta(\cal D)}\leq C\sup\limits_{z\in\cal D}|f(z)|.$$
\end{Lemma}

It is well known that every $m\times n$ $(m\leq n)$ matrix $A$ may be written
as $A=U(\sum\limits^m_{k=1}\lambda_k E_{kk})V,$ where
$U$ and $V$ are $m\times m$ and $n\times n$ unitary matrices respectively, and
$\lambda_1\geq\cdots\geq\lambda_m\geq 0,$ $E_{kk}$ is a $m\times n$ matrix,
the element of $k$th row and $k$th
column is 1, and other elements are $0.$ Hence
for every $P\in R_I(m,n)$ $ (m\leq n),$ there exist
$m\times m$ unitary matrix $U$ and $n\times n$ unitary matrix $V$,
such that $P=U(\sum\limits^m_{k=1}\lambda_k E_{kk})V,$
$(1\geq\lambda_1\geq\cdots\geq\lambda_m\geq 0).$

\begin{Lemma}([4], Lemma 8) Let $P=U\left(\begin{array}{ccccccc}\lambda_1&0&\cdots&0&0&\cdots&0\\
0&\lambda_2&\cdots&0&0&\cdots&0\\
\cdots&\cdots&\cdots&\cdots&\cdots&\cdots&\cdots\\
0&0&\cdots&\lambda_m&0&\cdots&0\end{array}\right)V\in R_I,$
and write $$Q=U\left(\begin{array}{cccc}\displaystyle\frac{1}
{\sqrt{1-\lambda_1^2}}&0&\cdots&0\\
0&\displaystyle\frac{1}
{\sqrt{1-\lambda_2^2}}&\cdots&0\\
\cdots&\cdots&\cdots&\cdots\\
0&0&\cdots&\displaystyle\frac{1}
{\sqrt{1-\lambda_m^2}}\end{array}\right)\overline{U}',$$
$$R=\overline{V}'\left(\begin{array}{ccccccc}\displaystyle\frac{1}
{\sqrt{1-\lambda_1^2}}&0&\cdots&0&0&\cdots&0\\
0&\displaystyle\frac{1}
{\sqrt{1-\lambda_2^2}}&\cdots&0&0&\cdots&0\\
\cdots&\cdots&\cdots&\cdots&\cdots&\cdots&\cdots\\
0&0&\cdots&\displaystyle\frac{1}
{\sqrt{1-\lambda_m^2}}&0&\cdots&0\\
&0_{(n-m)\times m}& & & &I_{n-m}&  \end{array}\right)V,$$ where
$U$ and $V$ are $m\times m$ and $n\times n$ unitary matrices respectively, and
$\lambda_1\geq\cdots\geq\lambda_m\geq 0.$ For $Z\in R_I,$ denote
$\Phi^{(I)}_P(Z)=Q(P-Z)(I_n-\overline{P}'Z)^{-1}R^{-1},$ then

(I) $\Phi^{(I)}_P\in Aut(R_I);$

(II) $\left(\Phi^{(I)}_P\right)^{-1}
=\Phi^{(I)}_P;$

(III) $\Phi^{(I)}_P(0)=0, \Phi^{(I)}_P(P)=P; $

(IV) $d\Phi^{(I)}_P(Z)|_{Z=P}=-QdZR,
d\Phi^{(I)}_P(Z)|_{Z=0}=-Q^{-1}dZR^{-1}; $

(V) $\Phi^{(I)}_P(Z)=Q^{-1}(I_m-Z\overline{P}')^{-1}(P-Z)R,$ for $Z\in R_I; $

(VI) $(I_m-Z\overline{P}')Q
\left(I_m-\Phi^{(I)}_P(Z)\overline{\Phi^{(I)}_P(Z)}'\right)\overline {Q}'
(I_m-P\overline{Z}')=I_m-Z\overline{Z}', $  for $Z\in R_I. $
\end{Lemma}

\section{An Important Proposition}
{\bf Proposition}\hspace*{2mm} Let $R_A (A=I,II,III,IV)$ be a classical
bounded symmetric domain. If $a^j\in R_A$, $d(a^j, \partial R_A)\to 0$ as
$j\to\infty,$ and
$w^j\in \mbox{\Bbb C}^N-\{0\},$ then exists a sequence of functions $\{f_{j}\}$
satisfying the following three conditions:

(i) $\{f_{j}\}$ is a bounded sequence in $\beta(R_A)$;

(ii) $\{f_{j}\}$ tends to zero uniformly on any compact subsets of $R_A;$

(iii) $\displaystyle\frac
{|\bigtriangledown f_j(a^j)w^j|}{H^{\frac{1}{2}}_{a^j}(w^j,w^j)}\geq C.$

{\bf Proof}\hspace{2mm} Note that the construction of the test
functions in [4], and replace $\phi(Z_j)$ by $a_j$ at a time, we
can construct a sequence of functions $\{f_{j}\}$ satisfying the
above three conditions. For example, for the reader's convenience,
we give the proof for the domain $R_I.$  But the proof completely
follows from [4], if necessary, the proof can be also omitted.

We construct a sequence of the functions according to the
following four parts.

Part A: To construct the sequence of $\{f_{j}\}$, we first assume
that
$$a^j=r_{j}E_{11},  j=1,2,\ldots , $$
where $E_{kl}$ is a $m\times n$ matrix, the element of $k$th row
and $l$th column is 1, and other elements are 0. It is clear that
$0<r_{j}<1$ and $r_{j}\to 1$ as $j\to\infty.$

Denote $w^{j}=(w^j_{11},
\cdots,w^j_{1n},w^j_{21},\cdots,w^j_{2n},\cdots,w^j_{m1},\cdots,w^j_{mn}).$
Using the formula (\ref{b}), we have
\begin{eqnarray*}&&H^{I}_{a^j}(w^{j}, w^{j})=
H^{I}_{r_jE_{11}}(w^{j}, w^{j})\\
&&=(m+n)w^j\left(\begin{array}{{lccr}}(1-r^2_j)^{-1}&0&\cdots&0\\
0&1&\cdots&0\\
\cdots&\cdots&\cdots&\cdots\\0&0&\cdots&1\end{array}\right)\times
\left(\begin{array}{{lccr}}(1-r^2_j)^{-1}&0&\cdots&0\\
0&1&\cdots&0\\
\cdots&\cdots&\cdots&\cdots\\0&0&\cdots&1\end{array}\right)\overline{w^j}'\\
&&=(m+n)\Biggl[\displaystyle\frac{|w^{j}_{11}|^{2}}{(1-r^2_j)^{2}}
+\displaystyle\frac{1}{1-r_{j}^2}\left(\sum\limits^{n}_{l=2}
|w^{j}_{1l}|^{2}+\sum\limits^{m}_{k=2}|w^{j}_{k1}|^{2}\right)
+\sum\limits_{2\leq k\leq m, 2\leq l\leq n}
|w^{j}_{kl}|^{2}\Biggr].\end{eqnarray*} Denote
\begin{eqnarray*}&&A^I_j=
\displaystyle\frac{|w^{j}_{11}|^{2}}{(1-r^2_j)^{2}}, \\
&&B^I_j=
\displaystyle\frac{1}{1-r_{j}^2}\left(\sum\limits^{n}_{l=2}
|w^{j}_{1l}|^{2}+\sum\limits^{m}_{k=2}|w^{j}_{k1}|^{2}\right),\\
&&C^I_j= \sum\limits_{2\leq k\leq m, 2\leq l\leq n}
|w^{j}_{kl}|^{2},\end{eqnarray*} then
\begin{equation}H^{I}_{a^j}(w^{j}, w^{j})=(m+n)(A^I_j+B^I_j+C^I_j).\label{7}
\end{equation}
We construct the functions according to three different cases:

Case 1\hspace{2mm}If for some $j,$
\begin{equation}\max(B^I_j, C^I_j)\leq A^I_j,\label{8}\end{equation}
then set
\begin{equation}f_{j}(Z)=\log\left(1-e^{-a(1-r_{j})}z_{11}\right)
-\log\left(1-z_{11}\right),\label{9}\end{equation} where
$Z=(z_{kl}), 1\leq k\leq m, 1\leq l\leq n$ and $a$ is any positive
number.

Case 2\hspace{2mm}If for some $j,$
\begin{equation}\max(A^I_j,C^I_j)\leq B^I_j,\label{10}\end{equation}
then set
\begin{equation}f_{j}(Z)=\left(\sum\limits^n_{l=2}
e^{-i\theta^{j}_{1l}}z_{1l}+\sum\limits^m_{k=2}e^{-i\theta^{j}_{k1}}z_{k1}
\right)\left(\displaystyle\frac{1}
{\sqrt{1-e^{-a(1-r_{j})}z_{11}}}-\displaystyle\frac{1}{\sqrt{1-z_{11}}}
\right),\label{11}\end{equation} where $a$ is any positive number,
and $\theta^{j}_{1l} =\arg w^{j}_{1l},$ $\theta^{j}_{k1} =\arg
w^{j}_{k1}.$ If $w^{j}_{1l}=0$ for some $l$ or $w^{j}_{k1}=0$ for
some $k$, replace the corresponding term
$e^{-i\theta^{j}_{1l}}z_{1l}$ or $e^{-i\theta^{j}_{k1}}z_{k1}$ by
$0.$

Case 3 \hspace{2mm}If for some $j,$
\begin{equation}\max(A^I_j,B^I_j)\leq C^I_j,\label{12}\end{equation}
then set
\begin{equation}f_{j}(Z)=\left(\sum\limits_{2\leq k\leq m, 2\leq l\leq n}e^{-i\theta^{j}_{kl}}
z_{kl}\right)\sqrt{1-z_{11}}\left(\displaystyle\frac{1}
{\sqrt{1-e^{-a(1-r_{j})}z_{11}}}-\displaystyle\frac{1}{\sqrt{1-z_{11}}}
\right),\label{13}\end{equation} where $a$ is any positive number,
and $\theta^{j}_{kl} =\arg w^{j}_{kl}, 2\leq k\leq m, 2\leq l\leq
n.$ If $w^{j}_{kl}=0$ for some $k$ or $l,$ replace the
corresponding term $e^{-i\theta^{j}_{kl}}z_{kl}$ by $0.$

Let $E$ be a compact subset of $R_I,$ then there exists a $\rho\in
(0,1)$ such that $|z_{11}|\leq\rho,$ for any $Z=(z_{kl})\in E.$ It
is easy to show that the sequence of functions defined by
(\ref{9}), (\ref{11}) and (\ref{13}) respectively, converges to
zero uniformly on $E$ as $j\to\infty,$ so the sequence satisfies
the conditions (ii).

Now we prove the above sequence satisfies the condition (i).

For the functions defined by (\ref{9}), it is easy to see
\begin{eqnarray*}\bigtriangledown f_j(Z)&=&
\left(\displaystyle\frac{\partial f_j}{\partial z_{11}}(Z),\cdots,
\displaystyle\frac{\partial f_j}{\partial z_{1n}}(Z),\cdots,
\displaystyle\frac{\partial f_j}{\partial z_{m1}}(Z),\cdots,
\displaystyle\frac{\partial f_j}{\partial z_{mn}}(Z)\right)\\
&=&\left(\displaystyle\frac{\partial f_j}{\partial z_{11}}(Z),
\cdots,0,\cdots,0,\cdots,0\cdots,0\right).\end{eqnarray*} From
formula (2), it is easy to know that the metric matrix of
$R_I(m,n)$ is
$$T(Z,Z)=(m+n)(I_{m}-Z{\overline Z}')^{-1}\times
(I_{n}-{\overline Z}'Z)^{-1},$$ so $T(0,0)=(m+n)I_{mn}, $ and
$$T^{-1}(Z,Z)=(m+n)^{-1}(I_{m}-Z{\overline Z}')\times
(I_{n}-{\overline Z}'Z).$$ Thus
\begin{eqnarray*}&&\bigtriangledown f_j(Z)T^{-1}(Z, Z)
\overline{\bigtriangledown f_j(Z)}'=(m+n)^{-1}
\left|\displaystyle\frac{\partial f_j}{\partial
z_{11}}(Z)\right|^2 \left(1-\sum\limits^n_{l=1}|z_{1l}|^2\right)
\left(1-\sum\limits^n_{k=1}|z_{k1}|^2\right)\\
&&=(m+n)^{-1}\left(1-\sum\limits^n_{l=1}|z_{1l}|^2\right)
\left(1-\sum\limits^n_{k=1}|z_{k1}|^2\right)\left|
\displaystyle\frac{-e^{-a(1-r_j)}}{1-e^{-a(1-r_j)}z_{11}}+
\displaystyle\frac{1}{1-z_{11}}\right|^2\\
&&\leq
(m+n)^{-1}(1-|z_{11}|^2)^2\left(\displaystyle\frac{2}{1-|z_{11}|}\right)^2
\leq 4(m+n)^{-1}(1+|z_{11}|)^2\leq 16(m+n)^{-1}.\end{eqnarray*}
Now Lemma 4 gives
$$\|f_j\|_{\beta(R_I)}\leq 16C(m+n)^{-1}.$$
This proves that the functions (\ref{9}) satisfy condition (i).

For the functions defined by (\ref{11}), If $Z\in R_{I}(m,n),$ we
have
$$I_{m}-Z{\overline Z}'
=\left(\delta_{st}-\sum\limits^{n}_{k=1}z_{sk}{\overline z}_{tk}
\right)_{1\leq s,t\leq m}>0.$$
$$I_{n}-Z'{\overline Z}
=\left(\delta_{st}-\sum\limits^{m}_{k=1}z_{ks}{\overline z}_{kt}
\right)_{1\leq s,t\leq n}>0.$$ Hence
\begin{equation}\delta_{11}-\sum\limits^{n}_{l=1}z_{1l}{\overline z}_{1l}
=1-\sum\limits^{n}_{l=1}|z_{1l}|^{2}>0.\label{14}\end{equation}
\begin{equation}\delta_{11}-\sum\limits^{m}_{k=1}z_{k1}\overline z_{k1}
=1-\sum\limits^{m}_{k=1}|z_{k1}|^{2}>0.\label{15}\end{equation}
Now (\ref{14}) and (\ref{15}) imply
\begin{eqnarray}|f_{j}(Z)|&=&\left|\left(\sum\limits^n_{l=2}
e^{-i\theta^{j}_{1l}}z_{1l}+\sum\limits^m_{k=2}e^{-i\theta^{j}_{k1}}
z_{k1}\right)\left(\displaystyle\frac{1}
{\sqrt{1-e^{-a(1-r_{j})}z_{11}}}-\displaystyle\frac{1}{\sqrt{1-z_{11}}}
\right)\right|\nonumber\\
&\leq&\left(\sum\limits^{n}_{l=2}|z_{1l}|+\sum\limits^{m}_{k=2}
|z_{k1}|\right)
\left(\left|\displaystyle\frac{1}{\sqrt{1-e^{-a(1-r_{j})}z_{11}}}\right|+
\left|\displaystyle\frac{1}{\sqrt{1-z_{11}}}\right|\right)\nonumber\\
&<&\left(\sqrt{n}\left(\sum\limits^n_{l=2}|z_{1l}|^2\right)^{\frac{1}{2}}
+\sqrt{m}\left(\sum\limits^m_{k=2}|z_{k1}|^2\right)^{\frac{1}{2}}
\right) \left(\displaystyle\frac{1}{\sqrt{1-|z_{11}|}}+
\displaystyle\frac{1}{\sqrt{1-|z_{11}|}}\right)\nonumber\\
&<&(\sqrt{n}+\sqrt{m})\sqrt{1-|z_{11}|^2}
\left(\displaystyle\frac{1}{\sqrt{1-|z_{11}|}}+
\displaystyle\frac{1}{\sqrt{1-|z_{11}|}}\right)\leq
4(\sqrt{m}+\sqrt{n}). \label{16}\end{eqnarray} From Lemma 5,
(\ref{16}) means that $\{f_{j}\}$ satisfy the condition (i).
Similarly, we may prove the functions defined by (\ref{13})
satisfy the condition (i).

At last, we prove that the sequence of functions defined by
(\ref{9}) satisfies the conditions (iii). In fact, by (\ref{7})
and (\ref{8}),
$$H^{I}_{a^j}(w^{j}, w^{j})=(m+n)(A^I_j+B^I_j+C^I_j)
\leq 3(m+n)A^I_j.$$
\begin{eqnarray*}&&\displaystyle\frac{\left|\bigtriangledown f_{j}(a^j)w^{j}
\right|}{H^{\frac{1}{2}}_{a^j}(w^j,w^j)}
\geq\sqrt{\displaystyle\frac{1}{3(m+n)}}
\displaystyle\frac{\left|\frac{\partial f_{j}} {\partial
z_{11}}(r_{j}E_{11})w^{j}_{11}\right|}
{\frac{|w^{j}_{11}|}{1-r^2_j}}\\
&&=\sqrt{\frac{1}{3(m+n)}}(1-r^2_j)
\left(\displaystyle\frac{1}{1-r_j}-\displaystyle\frac{e^{-a(1-r_{j})}}
{1-e^{-a(1-r_{j})}r_{j}}\right)\\
&&\geq\sqrt{\frac{1}{3(m+n)}}
\left(1-\displaystyle\frac{(1-r_{j})e^{-a(1-r_{j})}}
{1-e^{-a(1-r_{j})}r_{j}}\right),\end{eqnarray*} and
$$\lim\limits_{j\to\infty}\left[1-
\displaystyle\frac{(1-r_{j})e^{-a(1-r_{j})}}
{1-e^{-a(1-r_{j})}r_{j}}\right]=\displaystyle\frac{a}{a+1} \neq
0.$$ This proves the sequence of functions defined by (\ref{9})
satisfies the conditions (iii). Similarly, we can prove the
sequence of functions defined by (\ref{11}) or (\ref{13})
satisfies the conditions (iii).

Part B: We assume that $$a^j=r^{(1)}_jE_{11}+r^{(2)}_{j}E_{22},
$$ where $1>r^{(1)}_{j}\geq r^{(2)}_{j}\geq 0.$
By $a^j\to\partial R_{I},$ we may assume that $r^{(1)}_j\to 1,
r^{(2)}_{j}\to \lambda_0 (\leq 1).$

If $\lambda_0=1,$ using the same methods as in Part A, we can
construct a sequence of functions $\{f_{j}(Z)\}$ satisfying the
three conditions (i), (ii) and (iii).

If $\lambda_0<1,$ by Lemma 6, there exist
$\Phi^{I}_{r^{(1)}_jE_{11}+r^{(2)}_jE_{22}}\in R_I$ and
$\Phi^{I}_{r^{(1)}_jE_11}\in R_I,$ such that
$\Phi^{I}_{r^{(1)}_jE_{11}+r^{(2)}_jE_{22}}(r^{(1)}_jE_{11}+r^{(2)}_jE_{22})=0
$ and $\Phi^{I}_{r^{(1)}_jE_11}(r^{(1)}_jE_{11})=0$
($j=1,2,\cdots$). If we denote
$\Psi^{(j)}(Z)=\left(\Phi^{I}_{r^{(1)}_jE_{11}}\right)^{-1}
\circ\Phi^{I}_{r^{(1)}_jE_{11}+r^{(2)}_jE_{22}},$ then
$\Psi^{(j)}\in R_I$ and $\Psi^{(j)}\left(a^j\right)
=\Psi^{(j)}\left(r^{(1)}_jE_{11}+r^{(2)}_{j}E_{22}\right)=r^{(1)}_jE_{11}
=r_j E_{11},$ where $r_j=r^{(1)}_j.$

Set $g_{j}=f_{j}\circ\Psi^{(j)},$ where $\{f_j\}$ are the
functions obtained in Part A. Since $\Psi^{(j)}(Z)\in Aut(R_I),$
it is clear that
\begin{equation}H^I_{a^j}(w^j, w^j)=H^I_{\Psi^j(a^j)}
\left(J\Psi^{(j)}(a^j)w^{j}, J\Psi^{(j)}(a^j)w^{j}\right)
=H^I_{r_jE_{11}}(v^j, v^j),\label{18}\end{equation} where
$w^j=J\phi(Z^j)u^j, v^{j}=J\Psi^{(j)}(\phi(Z^{j}))w^j.$ It follows
from (\ref{18}) that
$$\displaystyle\frac{\left|\bigtriangledown (g_{j})
(a^j)w^{j}\right|}{H^{\frac{1}{2}}_{a^j}(w^{j}, w^{j})}
=\displaystyle\frac{\left|\bigtriangledown(f_{j})(r_{j}E_{11})
J\Psi^{(j)}(\phi(Z^{j}))w^{j})\right|}
{H^{\frac{1}{2}}_{r_{j}E_{11}}\left(
J\Psi^{(j)}(\phi(Z^{j}))w^{j},
J\Psi^{(j)}(\phi(Z^{j}))w^{j}\right)}=
\displaystyle\frac{\left|\bigtriangledown(f_{j})(r_{j}E_{11})
v^{j})\right|} {H^{\frac{1}{2}}_{r_{j}E_{11}}\left(v^{j},
v^{j}\right)}.$$ Now the discussion in Part A shows that
$$\displaystyle\frac{\left|\bigtriangledown(f_{j})(r_{j}E_{11})
v^{j})\right|} {H^{\frac{1}{2}}_{r_{j}E_{11}}\left(v^{j},
v^{j}\right)}\geq C>0,$$ that is, $\{g_{j}\}$ satisfies condition
(iii).

We prove that $\{g_{j}\}$ is a bounded sequence in $\beta(R_{I})$;
In fact, since $\Psi^{(j)}(Z)\in Aut(R_I),$
$$Q_{g_{j}}(Z)=Q_{f_{j}\circ\Psi^{(j)}}(Z)=Q_{f_{j}}\left(\Psi^{(j)}(Z)\right), $$
so $\|g_{j}\|_{\beta(R_{I})}=\|f_{j}\|_{\beta(R_{I})}$ is bounded.

Now we prove $\{g_{j}\}$ tends to zero uniformly on compact subset
$E$ of $R_{I}.$

If we write $\Psi^{(j)}(Z) =\left(\psi^{(j)}_{lk}(Z)\right)_{1\leq
l\leq m, 1\leq k\leq n}, $ by the definition of $\Psi^{(j)}$ and
Lemma 6, a direct calculation shows that
\begin{equation}\psi^{(j)}_{11}(Z)=z_{11}+r^{(2)}_j \displaystyle\frac{z_{12}z_{21}}
{1-r^{(2)}_{j}z_{22}}.\label{19}\end{equation} It is easy to show
that $\psi^{(j)}_{11}(Z)$ converges uniformly to
$\psi_{11}(Z)=z_{11}+\lambda_0 \displaystyle\frac{z_{12}z_{21}}
{1-\lambda_{0}z_{22}}$ on $R_{I}(m,n).$

Since $\lambda_0<1,$ $\lambda_0E_{11}+\lambda_0E_{22}\in R_I,$
similarly, there exists $\Psi(Z)\in Aut(R_{I}),$ such that
$\Psi(\lambda_0E_{11}+\lambda_0E_{22}) =\lambda_0E_{11},$ and the
first component of $\Psi(Z)$ is $\psi_{11}(Z).$ It is clear that
$\psi_{11}(Z)$ is holomorphic on $R_{I}.$ Let
$M_{1}=\sup\limits_{Z\in E} |\psi_{11}(Z)|=|\psi_{11}(Z_0)|,
(Z_0\in E).$ From $\Psi(Z)\in Aut(R_{I}),$ we know
$M_1=|\psi_{11}(Z_0)|<1,$ so we may choose $M_0>0$ with
$M_1<M_0<1.$ So for $j$ large enough,
$|\psi^{(j)}_{11}(Z_0)|<M_0,$ from this it follows that
$$1-|\psi^{(j)}_{11}(Z)|>1-M_0>0,$$
by the definition of $f_{j}(Z),$ it is easy to know
$g_{j}(Z)=f_{j}\circ\Psi^{(j)}(Z)$ tends to zero uniformly on $E.$
So $g_{j}(Z)$ satisfies the three conditions (i), (ii) and (iii).

Part C: Assume that
$$\phi(Z^j)=\sum\limits^m_{k=1}r^{(k)}_{j}E_{kk}, (1>r^{(1)}_j\geq
r^{(2)}_j\geq\ldots\geq r^{(m)}_j\geq 0). $$ By the same
discussion as in Part B, we can construct a sequence of functions
$f_{j}(z)$ which satisfies the three conditions (i), (ii) and
(iii).

Part D: In the general situation, $\phi(Z^j)\in R_{I}(m,n),$ so
there exist $m\times m$ unitary matrix $P_{j}$ and $n\times n$
unitary matrix $Q_{j}$, such that
$$P_{j}\left(\phi(Z^{j})\right)Q_{j}=
\sum\limits^m_{k=1}r^{(k)}_{j}E_{kk}. $$ We may assume that
$P_{j}\to P$ and $Q_{j}\to Q,$ as $j\to\infty$ (let
$P_j=(p^{kl}_{j}), P=(p^{kl}),$ $P_j\to P$ means that $p^{kl}_j\to
p^{kl}$ as $j\to\infty$ for any $1\geq k\geq m, 1\geq l\geq n$).
Let $\psi^{j}(Z)=P_jZQ_j,$ $\Psi(Z)=PZQ, Z\in R_{I}(m,n).$ It is
easy to show that $P$ is a $m\times m$ unitary matrix, $Q$ is a
$n\times n$ unitary matrix and $\psi^{(j)}(Z)$ converges uniformly
to $\psi(Z)$ on $R_{I}.$

Let $g_{j}(Z)=f_{j} \left(\psi^{(j)}(Z)\right),$ where $\{f_j\}$
are the functions obtained in Part C. From the same discussion as
that of Part B, we know $g_{j}(Z)$ satisfies conditions (i) and
(iii). For the compact subset $E\subset R_{I},$ it is easy to know
$\psi(E)$ is also a compact subset of $R_{I},$ so we can choose an
open subset $D_1$ of $R_I$ such that $\psi(E)\subset
D_1\subset\overline{D_1}\subset R_I.$ Since $\psi^{(j)}(Z)$
converges uniformly to $\psi(Z)$ on $R_{I}$, so as $j\to\infty,$
$\psi^j(E)\subset D_1.$
 Since $f_{j}(Z)$ tends to zero
uniformly on $\overline {D_1},$ we know $g_{j}(Z)=f_{j}
\left(\psi^{(j)}(Z)\right)$ tends to zero uniformly on $E\subset
R_I,$ i.e., $g_{j}$ satisfies condition (iii).

The last claim follows from the discussion in the above. This
completes the proof.

\section{Proof of Theorem}
It is well known that a bounded symmetric domain
except the two exceptional domains can be expressed as a
topological product of the first four types of irreducible domains
$R_A (A=I,II,III,IV)$. So we may assume that $\Omega=\Omega_{1}
\times\Omega_{2}\times\ldots\times\Omega_{n},$
where $\Omega_{k}\subset\mbox{\Bbb C}^{N_{k}}, (k=1,2\ldots,n)$ are
the classical bounded symmetric domains
and $N_{1}+N_{2}+\ldots+N_{n}=N.$ It is obvious $\Omega$ is homogeneous,
so by Lemma 3,
we only need to prove the condition (\ref{f}) is necessary.

Let $K_{\Omega}(z,z)$ and $T_{\Omega}(z,z)$ be the Bergman kernel function
and the Bergman metric of $\Omega$ respectively. Denote
$z=(z_{1}, z_{2},\ldots, z_{n})
\in\Omega\subset\mbox{\Bbb C}^{N},$ where $z_{k}\in\Omega_{k}\subset
\mbox{\Bbb C}^{N_{k}}.$
It is known that $$K_{\Omega}(z,z)=K_{\Omega_{1}}(z_{1}, z_{1})\times
K_{\Omega_{2}}(z_{2}, z_{2})\times \ldots
\times K_{\Omega_{n}}(z_{n}, z_{n}),$$
$$T_{\Omega}(z,z)=\left[\begin{array}{cccc}T_{\Omega_{1}}(z_{1}, z_{1})&
0&\cdots&0\\
0&T_{\Omega_{2}}(z_{2}, z_{2})&\cdots&0\\
\cdots&\cdots&\cdots&\cdots\\
0&0&\cdots&T_{\Omega_{n}}(z_{n}, z_{n})\end{array}\right]$$
\begin{equation}H^{\Omega}_{z}(u,u)=H^{\Omega_{1}}_{z_{1}}(u_{1},u_{1})+
H^{\Omega_{2}}_{z_{2}}(u_{2},u_{2})+\ldots+
H^{\Omega_{n}}_{z_{n}}(u_{n},u_{n}),\label{20}\end{equation}
where $u=(u_1, \ldots, u_n), u_k\in\mbox{\Bbb C}^k  (k=1,2, \ldots, n),$ and
$H^{\Omega}_z(u,u)=H_z(u,u).$

Now assume the condition (\ref{f}) fails, then there exists a sequence $\{z^{j}\}$
in $\Omega$ with $\phi(z^{j})\to\partial \Omega,$ as $j\to\infty,
u^{j}\in\mbox{\Bbb C}^{N}-\{0\},$ and an $\varepsilon_{0},$ such that
\begin{equation}\displaystyle\frac{H^{\Omega}_{\phi(z^{j})}
\left(J\phi(z^{j})u^{j}, J\phi(z^{j})u^{j}\right)}
{H^{\Omega}_{z^{j}}(u^{j},u^{j})}\geq\varepsilon_{0}, \label{21}\end{equation}
for all $j=1,2,\ldots$ . Write $J\phi(z^j)u^j=w^j,$
by (\ref{20}), we have
\begin{equation}H^{\Omega}_{\phi(z^j)}(w^j, w^j)
=H^{\Omega_{1}}_{\phi_1(z^j_{1})}(w^j_{1}, w^j_{1})+
H^{\Omega_{2}}_{\phi_2(z^j_{2})}(w^j_{2}, w^j_{2})+\ldots+
H^{\Omega_{n}}_{\phi_n(z^j_{n})}(w^j_{n}, w^j_{n}),\label{22}\end{equation}
where $w^j=(w^j_1, w^j_2,\ldots, w^j_n), w^j_k\in\mbox{\Bbb C}^{N_k},
 (k=1,2,\ldots,n),$ and $\phi=(\phi_1,\phi_2,\ldots,\phi_n).$
It is obvious that for some $k$, without loss of generality,
we may assume for $k=1,$ there exists a subsequence $\{j_s\}$ which we still
denote by $\{j\},$
such that
$$H^{\Omega_{l}}_{\phi_l(z^j_{l})}(w^j_{l}, w^j_{l})
\leq H^{\Omega_1}_{\phi_1(z^j_{1})}(w^j_{1}, w^j_{1}), \hspace*{4mm}
(l=2,\ldots, n).$$
So by (\ref{22})
\begin{equation}H^{\Omega}_{\phi(z^j)}(w^j, w^j)
\leq n H^{\Omega_{1}}_{\phi_1(z^j_{1})}(w^j_{1}, w^j_{1}) .
\label{23}\end{equation}

Since $\Omega_{1}$ is a classical bounded symmetric domain, denote
$a^j=\phi_1(z^j), $ by the important Proposition we may construct
a sequence of functions $\{f_{j}\}$
satisfying the following three conditions:

(a) $\{f_{j}\}$ is a bounded sequence in $\beta(\Omega_1)$;

(b) $\{f_{j}\}$ tends to zero uniformly on any compact subsets of $\Omega_1;$

(c) $\displaystyle\frac
{|\bigtriangledown f_j(a^j)w^j_1|}{H^{\frac{1}{2}}_{a^j}(w^j_1,w^j_1)}\geq C.$

We first prove that $\{f_{j}\}$ as a function sequence on $\Omega$
is a bounded sequence in $\beta(\Omega)$. In fact,
\begin{eqnarray*}\displaystyle\frac{\left|\bigtriangledown (f_{j})(z)u\right|}
{H^{\frac{1}{2}}_{z}(u, u)}&=&
\displaystyle\frac{\left|\bigtriangledown (f_{j})(z_1)u_1\right|}
{\left(H^{\Omega_1}_{z_1}(u_1, u_1)+\ldots+H^{\Omega_n}_{z_n}(u_n, u_n)\right)
^{\frac{1}{2}}}\\
&\leq&
\displaystyle\frac{\left|\bigtriangledown (f_{j})(z_1)u_1\right|}
{\left(H^{\Omega_1}_{z_1}(u_1, u_1)\right)^{\frac{1}{2}}}\leq \|f_j\|
_{\beta(\Omega_1)},\end{eqnarray*}
so $\{f_{j}\}$ is a bounded sequence in $\beta(\Omega)$.

For the compact subset $E\subset\Omega,$ let $E_1=\{z_1: (z_1, z_2, \ldots, z_n)
\in E\subset\Omega\}.$ It is easy to know $E_1$ is a compact subset of
$\Omega_1.$ We also have
$$\sup\limits_{z\in E}|f_j(z)|=\sup\limits_{z_1\in E_1}|f_j(z_1)|\to 0,$$ as
$j\to\infty.$ i.e., $\{f_{j}\}$ tends to zero uniformly on compact
subsets of $\Omega.$

Now we prove $\left\|C_{\phi}f_{j}\right\|_{\beta(\Omega)}\not\to 0.$
In fact, by (\ref{22}) and (\ref{23}), we have
\begin{eqnarray*}
&&\left\|C_{\phi}f_{j}\right\|_{\beta(\Omega)}=
\left\|f_{j}\circ\phi\right\|_{\beta(\Omega)} \geq Q_{f_{j}\circ\phi}(z^{j})\\
&&\geq \displaystyle\frac{\left|\bigtriangledown (f_{j}\circ\phi)(z^{j})u^{j}
\right|}{\left(H^{\Omega}_{z^{j}}(u,u)\right)^{\frac{1}{2}}}
=\displaystyle\frac{\left|\bigtriangledown (f_{j})(\phi(z^{j}))
J\phi(z^{j})u^{j}\right|}{\left(H^{\Omega}_{z^{j}}(u,u)\right)^{\frac{1}{2}}}\\
&&=\displaystyle\frac{\left|\bigtriangledown
(f_{j})(\phi(z^{j}))J\phi(z^{j})u^{j}\right|}
{\left(H^{\Omega}_{\phi(z^{j})}
\left(J\phi(z^{j})u^{j}, J\phi(z^{j})u^{j}\right)\right)^{\frac{1}{2}}}
\left\{\displaystyle\frac
{H^{\Omega}_{\phi(z^{j})}\left(J\phi(z^{j})u^{j}, J\phi(z^{j})u^{j}\right)}
{H^{\Omega}_{z^{j}}(u^{j},u^{j})}\right\}^{\frac{1}{2}}\\
&&\geq\sqrt{\varepsilon_{0}}
\displaystyle\frac
{\left|\bigtriangledown(f_{j})(\phi(z^j))w^{j}\right|}
{\left(H^{\Omega}_{\phi(z^j)}(w^{j},w^{j})\right)^{\frac{1}{2}}}
\geq\sqrt{\displaystyle\frac{\varepsilon_{0}}{n}}
\displaystyle\frac
{\left|\bigtriangledown(f_{j})(\phi_1(z^j))w^{j}_1\right|}
{\left(H^{\Omega_1}_{\phi_1(z^j)}(w^{j}_1, w^{j}_1)\right)^{\frac{1}{2}}}
=\sqrt{\displaystyle\frac{\varepsilon_{0}}{n}}\displaystyle\frac
{|\bigtriangledown f_j(a^j)w^j_1|}{H^{\frac{1}{2}}_{a^j}(w^j_1,w^j_1)}
,\end{eqnarray*}
by (c), we know $\left\|C_{\phi}f_{j}\right\|_{\beta(\Omega)}\not\to 0$ as
$j\to\infty.$
This contradicts the compactness of $C_{\phi}$ by Lemma 2.
Now the proof of the Theorem is complete.

\end{document}